\documentclass[11pt]{article}
\usepackage{amsthm,amsfonts,amsmath,amssymb}
\usepackage{mdwlist}
\newtheorem{thm}{Theorem}[section]

\newtheorem{lem}[thm]{Lemma}

\def\demo{\noindent{\bf Proof}\hskip10pt}

\def\qed{\hfill $\Box$}
\def\lg{\langle}
\def\rg{\rangle}
\def\rr#1{\item[{\rm (#1)}]}

\def\CD{\mathcal{CD}}
\def\L{\mathcal{L}}

\begin{document}
\title{Finite $p$-groups with at most $p^2+p$ subgroups not in Chermak-Delgado lattice  \thanks{This work was supported by NSFC (No. 11971280 \&11771258)}}
\author{Guojie Liu\ \ Haipeng Qu\ \ Lijian An\thanks{Corresponding author. e-mail: anlj@sxnu.edu.cn}\\
School of Mathematics and Computer Science,
Shanxi Normal University\\
Taiyuan, Shanxi 030032, P. R. China\\
 }

\maketitle

\begin{abstract}

The Chermak-Delgado lattice of a finite group $G$ is a self-dual sublattice of the subgroup lattice of $G$.
%Recent years, some scholars investigated finite groups whose Chermak-Delgado lattice is isomorphic to a given self-dual lattice.
In this paper, we determine finite $p$-groups with at most $p^2+p$ subgroups not in Chermak-Delgado lattice.

\medskip


\end{abstract}

\baselineskip=16pt

\section{Introduction}

 Suppose that $G$ is a finite group, and $H$ is a subgroup of $G$. The Chermak-Delgado measure of $H$ (in $G$) is denoted by $m_G(H)$, and defined as $m_G(H)=|H|\cdot |C_G(H)|.$
The maximal Chermak-Delgado measure of $G$ is denoted by $m^*(G)$, and defined as $$m^*(G)=\max\{ m_G(H)\mid H\le G\}.$$  Let $$\mathcal{CD}(G)=\{ H\mid m_G(H)=m^*(G)\}.$$ Then the set $\mathcal{CD}(G)$ forms a sublattice of $\L(G)$, the subgroup lattice of $G$, and is
called the Chermak-Delgado lattice of $G$.
It was first introduced by Chermak and Delgado \cite{CD}, and revisited by Isaacs \cite{I}.  In the last years, there has been a growing interest in understanding this lattice (see e.g. [1-3], [9-14], [16-18], [22-28]).
%\begin{thm}{\rm (\cite{I})}\label{th1} Given a finite group $G$, let $\mathcal{CD}(G)=\{ H\mid m_G(H)=m^*(G)\}.$ Then
%
%{\rm (1)} $\mathcal{CD}(G)$ is a lattice of subgroups of $G$.
%
%{\rm (2)} If $H,K\in\mathcal{CD}(G)$, then $\lg H,K\rg=HK$.
%
%{\rm (3)} If $H\in\mathcal{CD}(G)$, then $C_G(H)\in\mathcal{CD}(G)$ and $C_G(C_G(H))=H$.
%\end{thm}
%
%$\mathcal{CD}(G)$ is
%called the Chermak-Delgado lattice of G
%By Theorem \ref{th1},

Notice that a Chermak-Delgado lattice is  sublattice of subgroup lattice. It is natural to ask the question: How many subgroups  are not in the Chermak-Delgado lattices of finite groups. Some special cases of this question are proposed and solved. Fasol${\rm \check{a}}$ and T${\rm \check{a}}$rn${\rm \check{a}}$uceanu \cite{GF} classified groups with at most two subgroups not in Chermak-Delgado lattice. D. Burrell, W. Cocke and R. McCulloch \cite{DB} classified groups with at most four subgroups not in Chermak-Delgado lattice, and showed that the only non-nilpotent group with at most five subgroups not in the Chermak-Delgado lattice is $S_3$. By \cite[Lemma 5.1]{DB}, a non-abelian nilpotent group with five
subgroups not in the Chermak-Delgado lattice is a $2$-group or a $3$-group (Note that the dihedral group of order 8 is the only non-abelian nilpotent group with five groups not in the Chermak-Delgado lattice by the main theorem of this paper). In this paper, we focus on finite $p$-groups with few subgroups not in Chermak-Delgado lattice.

%Further information about the subject can refer \cite{ABQW}.
%
%
%types of Chermak-Delgado lattices are given.

% \begin{thm}{\rm (\cite{ABQW})}\label{th=abqw}
% If $\mathcal{L}$ is a Chermak-Delgado lattice of a finite $p$-group $G$ such that both $G/Z(G)$ and $G'$ are elementary abelian, then are $\mathcal{L}^+$ and $\mathcal{L}^{++}$, where $\mathcal{L}^+$ is a mixed $3$-string with center component isomorphic to $\mathcal{L}$ and the remaining components being $m$-diamonds {\rm (}a lattice with subgroups in the configuration of an $m$-dimensional cube{\rm )}, $\mathcal{L}^{++}$ is a mixed $3$-string with center component isomorphic to $\mathcal{L}$ and the remaining components being lattice isomorphic to $\mathcal{M}_{p+1}$.
%\end{thm}

Following D. Burrell, W. Cocke and R. McCulloch \cite{DB}, we use $\delta_{\mathcal{CD}}(G)$ to denote the number of subgroups of $G$ not in $\mathcal{CD}(G)$. That is, $$\delta_{\mathcal{CD}}(G)=|\mathcal{L}(G)|-|\mathcal{CD}(G)|.$$

Our main result is stated as follows.

\medskip

\noindent {\bf Main Theorem.}  Let $G$ be a finite $p$-group. Then $\delta_{\mathcal{CD}}(G)\leqslant p^{2}+p$ if and only if $G$ is one of the following groups:

$(1)$ $\mathrm{C}_{p^{k}}$, where $k=1,2,\ldots,p^{2}+p$; (In this case, $\delta_{\mathcal{CD}}(G)=k$.)

$(2)$ $\mathrm{C}_{p^t}\times\mathrm{C}_{p}$, where $t=1,2,\ldots,p-1$; (In this case, $\delta_{\mathcal{CD}}(G)=t(p+1)+1$.)

$(3)$ $\mathrm{Q}_{8}$; (In this case, $\delta_{\mathcal{CD}}(G)=1$.)

$(4)$ $\lg a,b\mid a^{p^k}=b^p=1,a^b=a^{1+p^{k-1}}\rg$, where $k=2,3,\ldots,p+1$. (In this case, $\delta_{\mathcal{CD}}(G)=(k-1)(p+1)$.)

\medskip

%\noindent {\bf Theorem B.} Let $G$ be a non-ablian $p$-group, where $p>2$. Then $\delta_{\mathcal{CD}}(G)\leqslant p^{2}+p$ if and only if
%\begin{center}
%$G\cong \mathrm{M}_{p}(k,1)$, where $k=2,3,\ldots,p+1$.
%\end{center}
%\medskip
%
%\noindent {\bf Theorem C.}
%Let $G$ be a non-ablian $2$-group. Then $\delta_{\mathcal{CD}}(G)\leqslant 6$ if and only if
%\begin{center}
%$G\cong \mathrm{D}_{8}$, $G\cong \mathrm{Q}_{8}$ or $G\cong \mathrm{M}_{2}(3,1)$.
%\end{center}

\medskip

\section{Preliminaries}
Let $G$ be a finite $p$-group. We use $s_{k}(G)$ to denote the number of subgroups of order $p^k$ of $G$. We use $\mathrm{exp}(G)$ and $d(G)$ to denote the exponent and minimal cardinality of generating set of $G$, respectively. We use $\mathrm{C}_{p^{n}}$ and $\mathrm{C}_{p}^{m}$ to denote the cyclic group of order $p^n$ and the direct product of $m$ copies of $\mathrm{C}_{p}$, respectively. We use $G_{i}$ to denote $i$-th member of the lower central series of $G$. We use $[G/Z(G)]$ to denote $\{H\mid Z(G)\leq H\leq G\}$.
Let $$r(G)=\mbox{max}\{k\mid E\leq G, E\cong \mathrm{C}_{p}^{k}\} ~\mbox{and}~r_{n}(G)=\mbox{max}\{k\mid E\unlhd G, E\cong \mathrm{C}_{p}^{k}\}.$$
Then $r(G)$ and $r_{n}(G)$ are called the rank and normal rank of $G$, respectively.
We define $\Omega_{1}(G)=\langle a\in G\mid a^{p}=1\rangle.$
Let $$\mathrm{M}_{p}(n,m):=\langle a,b\mid a^{p^{n}}=b^{p^{m}}=1, [a,b]=a^{p^{n-1}}\rangle~(n\geqslant 2,~m\geqslant1)~
\mbox{and}$$
$$\mathrm{M}_{p}(n,m,1):=\langle a,b\mid a^{p^{n}}=b^{p^{m}}=c^{p}=1, [a,b]=c,[c,a]=[c,b]=1\rangle~(n\geqslant m\geqslant1).$$

The following properties in Theorem \ref{basic} are basic and are often used in this paper. We will not point out when we use them.

\begin{thm}{\rm \cite[Theorem 2.1]{An}}\label{basic} Suppose that $G$ is a finite group and $H,K\in\CD(G)$.
\begin{itemize*}
  \rr{1} $\lg H,K\rg=HK$. Hence a Chermak-Delgado lattice is modular.
  \rr{2} $C_G(H\cap K)=C_G(H)C_G(K)$.
   \rr{3} $C_G(H)\in\mathcal{CD}(G)$ and $C_G(C_G(H))=H$. Hence a Chermak-Delgado lattice is self-dual, and $Z(G)\leq H$.
  \rr{4} Let $M$ be the maximal member of $\mathcal{CD}(G)$. Then $M$ is characteristic in $G$ and $\mathcal{CD}(M)=\mathcal{CD}(G)$.
  \rr{5} The minimal member of $\mathcal{CD}(G)$ is characteristic, abelian, and contains $Z(G)$.
\end{itemize*}
\end{thm}

%\begin{thm}{\rm \cite[Lemma 2.2]{BYW}}\label{basic2} Let $G$ be a group of order $p^n$. $s_{k}(G)$ denoted by the number of subgroups of order $p^k$ of $G$, where $k=0,1,\cdots,n$. Then
%\begin{center}
%$s_{k}(G)\equiv 1\ (\mbox{mod}~p)$ for $k=1,\cdots,n$.
%\end{center}
%\end{thm}

%\begin{defi}{\rm \cite[p.210]{MYX}}
%Let $G$ be a finite $p$-group.
%$$r(G)=\mbox{max}\{k\mid E\leq G, E\cong \mathrm{C}_{p}^{k}\}$$
%defined as the rank of $G$.
%$$r_{n}(G)=\mbox{max}\{k\mid E\unlhd G, E\cong \mathrm{C}_{p}^{k}\}$$
%defined as the normal rank of $G$.
%\end{defi}

%\begin{defi}{\rm \cite{MYX}}
%Let $G$ be a group of order $p^n$. $s_{k}(G)$ denoted by the number of subgroups of order $p^k$ of $G$, where $k=0,1,\cdots,n$.
%\end{defi}
\begin{lem}{\rm \cite[Theorem 2]{YC}}\label{qjg}
Let $G$ be a finite $p$-group. Then $\{C_{G}(H)\mid H\leq G\}=[G/Z(G)]$ if and only if $G'=\langle a\rangle$ is cyclic and either $p>2$, or $p=2$ and $[a,G]\leq\langle a^4\rangle$.

\end{lem}

\begin{lem}{\rm \cite[Lemma 4]{YC}}\label{qjg2}
If $G$ is as in Theorem \ref{qjg}, $Z=Z(G)$, and $Z\leq H\leq G$, then $|G/Z|=|H/Z|\cdot|C_{G}(H)/Z|$.
\end{lem}

\begin{lem}{\rm \cite[Theorem 5.1]{LH}}\label{lem4}
Let $G$ be a group of order $p^n$ where $p$ is an odd prime, and suppose that $r_{n}(G)=2$. Then one of the following holds:

$(1)$ $G$ is metacyclic;

$(2)$ $G=\langle a,x,y\mid a^{p^{n-2}}=x^p=y^p=1,[a,x]=y,[x,y]=a^{ip^{n-3}},[y,a]=1\rangle$ and $n\geqslant 4$. Moreover, $i=1$ or a fixed quadratic non-residue modulo $p$;

$(3)$ $G\cong \mathrm{M}_{p}(1,1,1)*\mathrm{C}_{p^{n-2}}$, $n\geqslant 3$;

$(4)$ $G$ is a $3$-group of maximal class.
\end{lem}

%\begin{lem}{\rm \cite[Theorem 2.2.9]{MYX}}\label{p3}
%Suppose that $G$ is a non-abelian group of order $p^3$. Then
%\begin{center}
%$G\cong \mathrm{Q}_{8}$, $G\cong \mathrm{M}_{p}(2,1)$ or $G\cong \mathrm{M}_{p}(1,1,1)$.
%\end{center}
%
%\end{lem}

\begin{lem}\label{lem7}{\rm \cite[Theorem 5.16]{YB}}\label{number}
Let a $p$-group $G$ be neither cyclic nor a $2$-group of maximal class,
$k\in \mathbb{N}$ and $p^{k}<p^{m}=|G|$.
Then $$s_{k}(G)\equiv 1+p\ (\mbox{mod}~p^2).$$
\end{lem}

\begin{lem}{\rm \cite[Theorem 1.2]{YB1}}\label{mp}
Let $G$ be a  quasi-regular metacyclic group of order $p^n$ and exponent $p^{e}<p^{n}$, $m\leqslant n$.

$(1)$ If $m\leqslant n-e$, then $s_{m}(G)=\frac{p^{m+1}-1}{p-1}$;

$(2)$ If $n-e<m\leqslant e$, then $s_{m}(G)=\frac{p^{n-e+1}-1}{p-1}$;

$(3)$ If $m>e$, then $s_{m}(G)=\frac{p^{n-m+1}-1}{p-1}$.

\end{lem}

\begin{lem}{\rm \cite[Proposition 3.13]{YB1}}\label{1+p}
Suppose that $G$ is a group of order $p^n>p^3$ satisfying $s_{k}(G)=1+p$ for some fixed $k\in\{2,\cdots,n-2\}$. Then one of the following holds:

$(1)$ $G$ is abelian of type $(p^n,p)$;

$(2)$ $G\cong \mathrm{M}_{p}(n-1,1)$;

$(3)$ $p=2,~k=2$, and $G=\langle a,b|a^{2^{n-2}}=1,n>4,a^{2^{n-3}}=b^{2}=1,a^{b}=a^{-1}\rangle$ is metacyclic.

\end{lem}

\begin{lem}\label{lem9}{\rm \cite[Corollary 1.7]{YB}}
Let $G$ be a group of maximal class of order $2^n$. Then $G$ is isomorphic to one of the following groups:

$(1)$ the dihedral group $\mathrm{D}_{2^{n}}$ where $n\geqslant 3$;

$(2)$ the generalized quaternion group $\mathrm{Q}_{2^{n}}$ where $n\geqslant 3$;

$(3)$ the semidihedral group $\mathrm{SD}_{2^{n}}$ where $n\geqslant 4$.

\end{lem}

\begin{lem}{\rm \cite[\S 9, Exercise 10]{YB}}\label{Mxg3}
Let $G$ be a group of maximal class of order $3^n$. Then $G_{1}$ is either abelian or minimal non-abelian, where $G_{1}=C_{G}(G_2/G_4)$.

\end{lem}

%\begin{lem}\label{mxg3}{\rm \cite[Theorem 4.3.15]{GJL}}
%Let $G$ be a group of maximal class of order $3^n$. Let $G_{1}=C_{G}(G_{n-2})$.
%
%$(1)$ If $n=3$, then $\CD(G)=[G/Z(G)]$;
%
%$(2)$ If $n=4$, then $\CD(G)=\{G_{1}\}$;
%
%$(3)$ If $n=5$, then $\CD(G)=\{G_{1}\}$ or $\CD(G)=[G_{1}/Z(G_{1})]\cup\{G,Z(G)\}$;
%
%$(3)$ If $n>5$, then $\CD(G)=\{G_{1}\}$ or $\CD(G)=[G_{1}/Z(G_{1})]$.
%
%\end{lem}

\begin{lem} \label{weiyi}{\rm \cite[Lemma 2.4]{BYW}}
Let $G$ be a group of order $p^n$.

$(1)$ If $1<m<n$ and $s_{m}(G)=1$, then G is cyclic;

$(2)$ If $s_{1}(G)=1$, then $G$ is cyclic or generalized quaternion group.
\end{lem}

\begin{lem}\label{nona2}{\rm \cite[\S 10, Remark 6]{YB}}
Let $G$ be a non-abelian group of order $2^4$. Then one of the following holds:
$$(1) \ G\cong \mathrm{Q}_{2^4};\ (2)~G\cong \mathrm{D}_{2^4}; ~(3)~G\cong \mathrm{SD}_{2^4}; ~(4)~G\cong \mathrm{M}_{2}(3,1); ~(5)~G\cong \mathrm{M}_{2}(2,2);$$ $$(6)G\cong \mathrm{M}_{2}(2,1,1);~(7)~G\cong \mathrm{D}_{8}\times\mathrm{C}_{2};~(8)~G\cong \mathrm{Q}_{8}\times\mathrm{C}_{2};~(9)~G\cong \mathrm{D}_{8}*\mathrm{C}_{4}.$$

\end{lem}

\begin{lem}\label{lem11}{\rm \cite[Theorem 2.2]{ZJ}}
Let $G$ be a non-abelian $2$-group with exactly three involutions. Then one of the following holds:

$(1)$ $G$ is a metacyclic $2$-group;

$(2)$ If $W$ is a maximal normal abelian non-cyclic subgroup of exponent $\leqslant4$ in $G$, then $W\cong \mathrm{C}_{4}\times\mathrm{C}_{2}$ or $W\cong \mathrm{C}_{4}\times\mathrm{C}_{4}$, $W=\Omega_{2}(C_{G}(W))$and $C_{G}(W)$ is metacyclic.

Suppose that $G$ has no normal subgroup isomophic to $\mathrm{C}_{4}\times\mathrm{C}_{4}$. Then $G$ has a normal subgroup $W\cong \mathrm{C}_{4}\times\mathrm{C}_{2}$, $C=C_{G}(W)$ is abelian of type $(2^n,2),~n\geqslant 2$, and $G/C$ is isomorphic to a proper non-trivial subgroup of $\mathrm{D}_{8}$.

\end{lem}

\section{Proof of  Main Theorem}
\begin{lem}\label{lem2}
Let $G$ be a finite $p$-group. Then the number of subgroups of order $p$ of $G$ in $\mathcal{CD}(G)$ is at most one.
\end{lem}

\demo
Let $N\leq G$ with $|N|=p$. If $N\in \mathcal{CD}(G)$, then it is easy to see that $N=Z(G)$.\qed

\begin{lem}\label{lem3}
Let $G$ be a finite $p$-group. If $\delta_{\mathcal{CD}}(G)\leqslant p^2+p$, then $r(G)\leqslant 2$ and $G$ has no subgroup isomorphic to ${\rm M}_{p}(1,1,1)$.
\end{lem}

\demo
Otherwise, there exists $H\leq G$ such that $H$ is isomorphic to $\mathrm{C}_{p}^{3}$ or ${\rm M}_{p}(1,1,1)$. By calculation, $s_{1}(G)=\frac{p^{3}-1}{p-1}=p^{2}+p+1$. By Lemma \ref{lem2}, there are at least $p^2+p$ subgroups of order $p$ not in $\CD(G)$. Note that $1\notin \mathcal{CD}(G)$. We have $\delta_{\mathcal{CD}}(G)\geqslant p^2+p+1$. This contradicts the hypothesis.
\qed
%\mbox{\hskip .2 in}~Let $G$ be an abelian $p$-group. Then $\mathcal{CD}(G)=\{G\}$. It is easy to see that $\delta_{\mathcal{CD}}(G)=|\mathcal{L}(G)|-1$. That is

\begin{lem}\label{lem1}
Let $G$ be a finite abelian $p$-group. Suppose that $\delta_{\mathcal{CD}}(G)\leqslant p^2+p$. Then $G$ is one of the following groups:

 $(1)$ $\mathrm{C}_{p^{k}}$, where $k=1,2,\ldots,p^{2}+p$; $($In this case, $\delta_{\mathcal{CD}}(G)=k$.$)$

$(2)$ $\mathrm{C}_{p^s}\times\mathrm{C}_{p}$, where $s=1,2,\ldots,p-1$. $($In this case, $\delta_{\mathcal{CD}}(G)=s(p+1)+1$.$)$

\end{lem}
\demo
Since $G$ is abelian, $\mathcal{CD}(G)=\{G\}$. By Lemma \ref{lem3}, $r(G)\leqslant 2$.

If $r(G)=1$, then $G\cong \mathrm{C}_{p^{k}}$. In this case, $\delta_{\mathcal{CD}}(G)=|\mathcal{L}(G)|-|\CD(G)|=k+1-1=k$. Hence $1\leqslant k\leqslant p^{2}+p$.

If $r(G)=2$, then we may assume that $G\cong \mathrm{C}_{p^s}\times\mathrm{C}_{p^t}$ where $s\ge t$. We assert that $t$ is 1. If not, then there exists $K\leq G$ such that $K\cong \mathrm{C}_{p^2}\times\mathrm{C}_{p^2}$. By Lemma \ref{mp},
$$|\mathcal{L}(K)|=\sum_{i=0}^{4} s_i(K)=1+1+p+1+p+p^{2}+1+p+1=p^2+3p+5.$$
Hence $\delta_{\mathcal{CD}}(G)=|\mathcal{L}(G)|-|\CD(G)|=p^2+3p+1$. This contradicts the hypothesis. Thus $G\cong \mathrm{C}_{p^s}\times\mathrm{C}_{p}$.
By Lemma \ref{mp},
$$|\mathcal{L}(G)|=\sum_{i=0}^{s+1} s_i(G)=1+s(1+p)+1=s(1+p)+2.$$
In this case, $\delta_{\mathcal{CD}}(G)=|\mathcal{L}(G)|-|\CD(G)|=s(p+1)+1$. Hence $1\leqslant s<p$.
\qed

%\demo
%By Lemma \ref{mx3}, $G_{1}$ is abelian or minimal non-abelian. If $G_{1}$ is abelian, then $G$ has a abelian maximal subgroup. If $n=3$, then $\CD(G)=[G/Z(G)]$. If $n>3$, then $|G/Z(G)|\geqslant p^{3}$. Hence $\CD(G)=\{G_{1}\}$.
%
%If $G_{1}$ is minimal non-abelian, then $n\geqslant5$ and $\CD(G_{1})=[G_{1}/Z(G_{1})]$.
%由定理\ref{thm5.3.4}可得,
%$\CD(G)=\CD(G_{1})$\\$\cup\{G,Z(G)\}$或$\CD(G)=\CD(G_{1})$.
%
%当$n=5$时, $m_{G}(G)=|G||Z(G)|=p^{6}=|G_{1}||Z(G_{1})|=|G_{1}||C_{G}(G_{1})|=m_{G}(G_{1})$. 则$G\in \CD(G)$. 故$\CD(G)=[G_{1}/Z(G_{1})]\cup\{G,Z(G)\}$.
%
%当$n>5$时, $m_{G}(G)=|G||Z(G)|=3^{n+1}<3^{2n-4}=|G_{1}||Z(G_{1})|=|G_{1}||C_{G}(G_{1})|=m_{G}(G_{1})$. 故$G\notin \CD(G)$, 则$\CD(G)=[G_{1}/Z(G_{1})]$.\qed

\begin{thm}\label{yxh}
Let $G$ is metacyclic $p$-group, where $p>2$. Then $\CD(G)=[G/Z(G)]$.
\end{thm}
\demo
Since $G$ is metacyclic, $G'$ is cyclic.  By Lemma \ref{qjg}, $\{C_{G}(H)\mid H\leq G\}=[G/Z(G)]$. By Lemma \ref{qjg2}, $m_{G}(H)=|H|\cdot|C_{G}(H)|=|G|\cdot|Z(G)|=m_{G}(G)$ for any $H\in [G/Z(G)]$. Hence $\CD(G)=[G/Z(G)]$.

\begin{thm}\label{jdl3}
Let $G$ be a $3$-group of maximal class of order $3^n$. Let $G_{1}=C_{G}(G_{2}/G_{4})$.

$(1)$ If $n=3$, then $\CD(G)=[G/Z(G)]$ $($In this case, $|\CD(G)|=p+3$$)$;

$(2)$ If $n=4$, then $\CD(G)=\{G_{1}\}$;

$(3)$ If $n=5$, then $\CD(G)=\{G_{1}\}$ or $\CD(G)=[G_{1}/Z(G_{1})]\cup\{G,Z(G)\}$ $($In this case, $|\CD(G)|=p+5$$)$;

$(4)$ If $n>5$, then $\CD(G)=\{G_{1}\}$ or $\CD(G)=[G_{1}/Z(G_{1})]$ $($In this case, $|\CD(G)|=p+3$$)$.
\end{thm}
\demo
If $n=3$, then $G/Z(G)\cong \mathrm{C}_{p}\times \mathrm{C}_{p}.$ Thus $$\mathcal{CD}(G)=[G/Z(G)] ~\mbox{and}~ |\mathcal{CD}(G)|=|\mathcal{L}(G/Z(G))|=1+1+p+1=p+3.$$ By Lemma \ref{Mxg3}, $G_{1}$ is abelian or minimal non-abelian. If $G_{1}$ is abelian, then $G$ has an abelian maximal subgroup. If $n\geqslant 4$,
then $\mathcal{CD}(G)=\{G_{1}\}.$  If $G_{1}$ is minimal non-abelian, then $$n\geqslant 5~\mbox{and}~G_{1}/Z(G_{1})\cong \mathrm{C}_{p}\times \mathrm{C}_{p}.$$ Thus $$\mathcal{CD}(G_{1})=[G_{1}/Z(G_{1})] ~\mbox{and}~ |\mathcal{CD}(G_{1})|=|\mathcal{L}(G_{1}/Z(G_{1}))|=p+3.$$ For any $M\lessdot G$, $C_{G}(M)\trianglelefteq G$. Thus $C_{G}(M)\leqslant M$. Hence $C_{G}(M)=Z(M)$. By \cite[\S 9, Theorem 9.6(e)]{YB}, all maximal subgroups except $G_{1}$ are $3$-groups of maximal class. Therefore, for any $M\lessdot G$ and $M\neq G_{1}$, $$m_{G}(M)=|M|\cdot|Z(M)|=p^{n}<p^{2n-3}=|G_{1}|\cdot|Z(G_{1})|=m_{G}(G_{1}).$$ Hence $G_{1}$ is the unique maximal subgroup in $\CD(G)$. Thus $\mathcal{CD}(G_{1})\subseteq \mathcal{CD}(G)$. If $n=5$, then $m_{G}(G)=m_{G}(G_{1})$. Hence $$\CD(G)=[G_{1}/Z(G_{1})]\cup\{G,Z(G)\} ~\mbox{and}~ |\mathcal{CD}(G)|=p+5.$$
If $n>5$, then $m_{G}(G)<m_{G}(G_{1})$. Hence $\CD(G)=[G_{1}/Z(G_{1})]$.

\begin{lem}\label{zhu2}
Let $G$ be a finite non-ablian $p$-group, where $p>2$. Suppose that $\delta_{\mathcal{CD}}(G)\leqslant p^{2}+p$. Then $G\cong \mathrm{M}_{p}(k,1), \mbox{where} ~k=2,3,\ldots,p+1$. Moreover, $ \delta_{\mathcal{CD}}(G)=(k-1)(p+1).$
\end{lem}

\demo
By Lemma \ref{lem3}, $r(G)=2$ and $G$ has no subgroup isomorphic to $\mathrm{M}_{p}(1,1,1)$. Since $r(G)=2$, $G$ is one of the groups in Lemma \ref{lem4}.
If $G$ is of type (2) in Lemma \ref{lem4}, then $\langle a^{p^3},x,y\rangle\cong \mathrm{M}_{p}(1,1,1)$. Hence $G$ is one of type (1) or (4) in Lemma \ref{lem4}.

We assert that $G$ is metacyclic. Otherwise, $G$ is a $3$-group of maximal class of order $3^n$, where $n\geqslant 4$. By Theorem \ref{jdl3}, if $n=4$, then $|\mathcal{CD}(G)|=1$; if $n=5$, then $|\mathcal{CD}(G)|\leqslant p+5=8$; if $n>5$, $|\mathcal{CD}(G)|\leqslant p+3=6$. Since $G$ is not the group in Lemma \ref{1+p}, $s_{t}(G)\geqslant 1+2p,~t=2,\ldots,n-2$. Then if $n=4$, $|\mathcal{L}(G)|\geqslant 4p+5=17$; if $n=5$, $|\mathcal{L}(G)|\geqslant 6p+6=24$;
if $n>5$, $|\mathcal{L}(G)|\geqslant 8p+7=31$. Thus if $n=4$, $\delta_{\mathcal{CD}}(G)\geqslant 16>12$; if $n=5$, $\delta_{\mathcal{CD}}(G)\geqslant 16>12$; if $n>5$, $\delta_{\mathcal{CD}}(G)\geqslant 25>12$. This contradicts the hypothesis.

Hence $G$ is metacyclic. By Lemma \ref{yxh}, $\mathcal{CD}(G)=[G/Z(G)]$. Metacyclic $p$-groups, $p>2$, are quasi-regular (\cite[Theorem 7.2(c)]{YB}).

We assert that $n-e=1$. If not, by Lemma \ref{mp}, $s_{2}(G)=1+p+p^2$. (1) If $|Z(G)|\geqslant p^2$, then there exists $p^{2}+p$ subgroups of order $p^2$ and $p+1$ subgroups of order $p$ not in $\mathcal{CD}(G)$. In this case, $\delta_{\mathcal{CD}}(G)\geqslant p^2+2p+1$. This contradicts the hypothesis.
(2) If $|Z(G)|=p$, then $m^*(G)=p^{n+1}$. Hence there exist $p^{2}$ subgroups of order $p^2$ not in $\mathcal{CD}(G)$, since $G$ has only $1+p$ maximal subgroups and $\mathcal{CD}(G)$ is self-dual. By lemma \ref{lem2}, there exists $p$ subgroups of order $p$ not in $\mathcal{CD}(G)$. Note that $1\notin \mathcal{CD}(G)$. We have $\delta_{\mathcal{CD}}(G)\geqslant p^2+p+1$. This contradicts the hypothesis.

Hence $G$ has cyclic maximal subgroup. By Lemma \cite[Theorem 1.2]{YB}, $G\cong \mathrm{M}_{p}(k,1)$. Since $|G:Z(G)|=p^2$ and $\mathcal{CD}(G)=[G/Z(G)]$, $|\mathcal{CD}(G)|=p+3$. By Lemma \ref{mp},
$$|\mathcal{L}(G)|=\sum_{i=0}^{k+1} s_i(G)=1+k(1+p)+1=k(1+p)+2.$$ Thus $\delta_{\mathcal{CD}}(G)=(k-1)(p+1)$. Since $\delta_{\mathcal{CD}}(G)\leqslant p^{2}+p$, $2\le k\le p+1$.\qed

\begin{lem}\label{lem6}
Let $G$ be a non-abelian group of order $2^n$. If $\delta_{\mathcal{CD}}(G)\leqslant 6$, then $G\in \mathcal{CD}(G)$.

\end{lem}

\demo
Otherwise, let $M$ be the the maximal member of $\mathcal{CD}(G)$. Then $|M|\leqslant 2^{n-1}$ and $|M|\cdot|Z(M)|=m_G(M)>2^{n+1}$. Hence $|Z(M)|>4$. By Lemma \ref{weiyi}, there are at least three subgroups of order $4$ in $\L(G)\setminus\CD(G)$. Note that there exist at least two maximal subgroups and one subgroup of order $2$ not in $\mathcal{CD}(G)$, and $G,1\notin \mathcal{CD}(G)$.
We have $\delta_{\mathcal{CD}}(G)\geqslant 8$, which contradicts the hypothesis. \qed

\begin{lem}\label{cc}
Let $G$ be a non-abelian group of order $2^n$. If $\delta_{\mathcal{CD}}(G)\leqslant 6$, then $Z(G)\cong \mathrm{C}_{2}$ or $\mathrm{C}_4$. Moreover, if $Z(G)\cong \mathrm{C}_4$, then $|G|\leqslant 2^6$.
\end{lem}

\demo
We assert that $Z(G)$ is cyclic. If not, then there exists $N\leq Z(G)$ such that $N\cong \mathrm{C}_{2}\times\mathrm{C}_{2}$. It follows that $N$ has three subgroups of order $2$ not in $\mathcal{CD}(G)$.
Since $G$ is non-abelian, $G/N$ is not cyclic. Hence there exists $M\geq N$ such that $G/M\cong \mathrm{C}_{2}\times\mathrm{C}_{2}$. Let $\lg a_i,M\rg\ (i=1,2,3)$ be three maximal subgroup of $G$. Then $\langle a_{1}\rangle$, $\langle a_{2}\rangle$ and $\langle a_{3}\rangle$ are not in $\mathcal{CD}(G)$, since they do not contain $Z(G)$.
  Note that $1\notin \mathcal{CD}(G)$. We have $\delta_{CD}(G)\geqslant 7$, which contradicts the hypothesis. Hence $Z(G)\cong \mathrm{C}_{2^{k}}$.

If $k\geqslant2$, then, by Lemma \ref{number}, there exist at least three subgroups of order $4$, three subgroups of order $2$, since $G$ is not of maximal class. Note that $1$ not in $\mathcal{CD}(G)$. We have $\delta_{\mathcal{CD}}(G)\geqslant 6$, where ``=" holds if and only if $k=2$. Since $\delta_{\mathcal{CD}}(G)\leqslant 6$,
\begin{center}
$\delta_{\mathcal{CD}}(G)=6$ and $Z(G)\cong \mathrm{C}_{4}$.
 \end{center}
 In this case, $\Omega_1(G)\cong\mathrm{C}_2\times \mathrm{C}_2$. Let $K=Z(G)\Omega_1(G)\cong\mathrm{C}_4\times \mathrm{C}_2$. Then $K$ has three subgroups of order $4$ and three subgroups of order $2$. It forces that $K=\Omega_2(G)$. Let $\overline{G}=G/\Omega_{1}(G)$. Since $\Omega_1(\overline{G})=\Omega_2(G)/\Omega_1(G)\cong \mathrm{C}_2$, $\overline{G}$ has exactly one element of order $2$. By Lemma \ref{weiyi},
 \begin{center}
  $\overline{G}\cong \mathrm{C}_{2^{n-2}}$ or $\overline{G}\cong \mathrm{Q}_{2^{n-2}}$.
\end{center}

 If $\overline{G}\cong \mathrm{C}_{2^{n-2}}$, then there exists a cyclic subgroup of order $2^{n-1}$ of $G$. By Lemma \ref{lem6}, $G\in \mathcal{CD}(G)$. Then $m^{*}(G)=|G|\cdot|Z(G)|=2^{n+2}\geqslant 2^{n-1}\cdot 2^{n-1}$. Hence $n\leqslant 4$.

If $\overline{G}\cong \mathrm{Q}_{2^{n-2}}$, then there exists a cyclic subgroup of order $2^{n-2}$ of $G$. By Lemma \ref{lem6}, $G\in \mathcal{CD}(G)$. Then $m^{*}(G)=|G|\cdot|Z(G)|=2^{n+2}\geqslant 2^{n-2}\cdot 2^{n-2}$. Hence $n\leqslant 6$.
\qed

%\begin{lem}\label{lem8}
%Let $G$ be a non-cyclic group of order $2^n$ and not be of maximal class. If $\delta_{\mathcal{CD}}(G)\leqslant 6$, then $s_{1}(G)=3$.
%
%\end{lem}
%\demo
%If not, then by Lemma \ref{lem7}, $s_{1}(G)\geqslant 7$. By Lemma \ref{lem2} and $1\notin \mathcal{CD}(G)$, $\delta_{\mathcal{CD}}(G)\geqslant 7$, this is a contradiction.\qed

%\begin{lem}\label{m2}
%Let $G$ be a $2$-group of maximal class. If $\delta_{\mathcal{CD}}(G)\leqslant 6$, then $G\cong \mathrm{Q}_{8}$ or $G\cong \mathrm{D}_{8}$.
%
%\end{lem}
%\demo
%By Lemma \ref{lem9}, $G$ has a cyclic maximal subgroup. Thus if $n\geqslant 4$, then $G\notin \mathcal{CD}(G)$. This contradicts that Lemma \ref{lem6}. By calculation, $\delta_{\mathcal{CD}}(\mathrm{Q}_{8})=1$ and $\delta_{\mathcal{CD}}(\mathrm{D}_{8})=5$. Hence $G\cong \mathrm{Q}_{8}$ or $G\cong \mathrm{D}_{8}$.\qed

\begin{lem}\label{lem12}
Let $G$ be a non-abelian group of order $2^n$. If $\delta_{\mathcal{CD}}(G)\leqslant 6$ and $Z(G)\cong \mathrm{C}_{2}$, then $G$ has no normal subgroup isomorphic to $\mathrm{C}_{4}\times\mathrm{C}_{4}$.
\end{lem}

\demo
If not, let $H\unlhd G$ such that $H\cong\mathrm{C}_{4}\times\mathrm{C}_{4}=\langle a\rangle\times\langle b\rangle$.
By calculation, $H$ has six cyclic subgroups of order $4$. Since $H\unlhd G$, $Z(G)\leq H$. Let $Z(G)=\langle a^2\rangle$. By calculation, $\langle a\rangle\geq Z(G)$, $\langle ab^2\rangle\geq Z(G)$. That is $H$ has four cyclic subgroups of order $4$ not in $\mathcal{CD}(G)$. Hence $G$ has at least four cyclic subgroups of order $4$ not in $\mathcal{CD}(G)$. It follows from $G$ has at least two subgroups of order $2$ and one not in $\mathcal{CD}(G)$ that $\delta_{\mathcal{CD}}(G)\geqslant 7$. This contradicts the hypothesis.\qed

\begin{lem}\label{lemn7}
Let $G$ be a non-metacyclic group of order $2^n$. If $\delta_{\mathcal{CD}}(G)\leqslant 6$ and $Z(G)\cong \mathrm{C}_{2}$, then $|G|\leqslant 2^7$.
\end{lem}

\demo
If $G$ is $2$-group of maximal class, by Lemma \ref{lem9}, $G$ has a cyclic maximal subgroup. Thus if $n\geqslant 4$, then $G\notin \mathcal{CD}(G)$. This contradicts that Lemma \ref{lem6}. Hence $n\leqslant3$.

If $G$ is not of maximal class, we assert that $s_{1}(G)=3$. If not, then by Lemma \ref{lem7}, $s_{1}(G)\geqslant 7$. By Lemma \ref{lem2} and $1\notin \mathcal{CD}(G)$, $\delta_{\mathcal{CD}}(G)\geqslant 7$, This contradicts the hypothesis. Hence by Lemma \ref{lem11} and Lemma \ref{lem12}, there exists $N\unlhd G$ such that $N\cong \mathrm{C}_{4}\times\mathrm{C}_{2}$, $C_{G}(N)\cong \mathrm{C}_{2^k}\times\mathrm{C}_{2}(k\geqslant 2)$ and $|G:C_{G}(N)|\leqslant 2^3$. By Lemma \ref{lem6}, $G\in \mathcal{CD}(G)$. It follows that $2^{n+1}=|G|\cdot|Z(G)|=m^{*}(G)\geqslant m_{G}(C_{G}(N))\geqslant |C_{G}(N)|^2\geqslant 2^{2(n-3)}$. Hence $n\leqslant 7$.

As above, $n\leqslant 7$.\qed

\begin{lem}\label{lem13}
Let $G$ be a metacyclic group of order $2^n$. If $\delta_{\mathcal{CD}}(G)\leqslant 6$ and $Z(G)\cong \mathrm{C}_{2}$, then $|G|\leqslant 2^5$.
\end{lem}

\demo
Let $G=\langle a\rangle\langle b\rangle$, where $\langle a\rangle\unlhd G$. By Lemma \ref{lem6},  $G\in \mathcal{CD}(G)$. Then $$2^{n+1}=|G|\cdot|Z(G)|=m^{*}(G)\geqslant m_{G}(\langle a\rangle)\geqslant o(a)^2$$ and $$2^{n+1}=|G|\cdot|Z(G)|=m^{*}(G)\geqslant m_{G}(\langle b\rangle)\geqslant o(b)^2.$$ Hence $o(a)\leqslant 2^{\frac{n+1}{2}}$ and $o(b)\leqslant 2^{\frac{n+1}{2}}$. It follows that $$2^n=|G|=|\langle a\rangle\langle b\rangle|=\frac{|\langle a\rangle|\cdot|\langle b\rangle|}{|\langle a\rangle\cap\langle b\rangle|}\leqslant\frac{2^{\frac{n+1}{2}}\cdot2^{\frac{n+1}{2}}}{|\langle a\rangle\cap\langle b\rangle|}=\frac{2^{n+1}}{|\langle a\rangle\cap\langle b\rangle|}.$$ Thus, $|\langle a\rangle\cap\langle b\rangle|\leqslant 2$.

If $|\langle a\rangle\cap\langle b\rangle|=2$, then $\langle a\rangle\cap\langle b\rangle=Z(G)$, since $Z(G)\cong \mathrm{C}_{2}$ and $\langle a\rangle\unlhd G$. Notice that $o(a)=o(b)=2^{k}$, where $k=\frac{n+1}{2}$. By $\mathrm{Aut}(\mathrm{C}_{2^k})\cong \mathrm{C}_{2^{k-2}}\times\mathrm{C}_{2}$, $a^{b^{2^{k-2}}}=a$. Then $b^{2^{k-2}}\in Z(G)$. This contradicts that $Z(G)=\langle a^{2^{k-1}}\rangle=\langle b^{2^{k-1}}\rangle$.

Thus $|\langle a\rangle\cap\langle b\rangle|=1$. We assert that $o(b)\leqslant 4$. If not, then there exist at least two subgroups of order $2$, two subgroups of order $4$ and two subgroups of order $8$ that do not contain $Z(G)$. Notice that $1\notin\mathcal{CD}(G)$. Then $\delta_{\mathcal{CD}}(G)\geqslant 7$. This contradicts the hypothesis. Since $2^n=|G|=|\langle a\rangle|\cdot|\langle b\rangle|\leqslant o(a)\cdot 4$, $o(a)\geqslant 2^{n-2}$. It follows from $$2^{n+1}=|G|\cdot|Z(G)|=m^{*}(G)\geqslant m_{G}(\langle a\rangle)\geqslant o(a)^2\geqslant 2^{2(n-2)}$$ that $2n-4\leqslant n+1$. Thus $n\leqslant 5$. \qed

\begin{lem}\label{zhu3}
Let $G$ be a non-ablian $2$-group. If $\delta_{\mathcal{CD}}(G)\leqslant 6$, then $G$ is one of the following groups:

$(1)$ $\mathrm{Q}_{8}$ $(\mbox{In this case}~\delta_{\mathcal{CD}}(\mathrm{Q}_{8})=1)$;

$(2)$ $\mathrm{D}_{8}$ $(\mbox{In this case}~\delta_{\mathcal{CD}}(\mathrm{D}_{8})=5)$;

$(3)$ $\mathrm{M}_{2}(3,1)$ $(\mbox{In this case}~\delta_{\mathcal{CD}}(\mathrm{M}_{2}(3,1))=6)$.

\end{lem}

\demo
By Lemma \ref{cc}, Lemma \ref{lemn7} and Lemma \ref{lem13}, we just need to consider groups of order at most $2^7$.

%Case $1$: $|G|=2^3$. Then $G\cong \mathrm{Q}_{8}$ or $G\cong \mathrm{D}_{8}$. By calculation,  $\delta_{\mathcal{CD}}(\mathrm{Q}_{8})=1$ and $\delta_{\mathcal{CD}}(\mathrm{D}_{8})=5$.
%
%Case $2$: $|G|=2^4$. We assert that $G\cong \mathrm{M}_{2}(3,1)$. If not, by Lemma \ref{nona2}, $G$ is one of the groups in Lemma \ref{nona2} except for $\mathrm{M}_{2}(3,1)$. If $G\cong \mathrm{Q}_{2^4}, ~G\cong \mathrm{D}_{2^4}~ \mbox{or}~ G\cong \mathrm{SD}_{2^4}$, then there exists a cyclic maximal subgroups of $G$ and $|Z(G)|=2$. Thus $G\notin \mathcal{CD}(G)$, this contradicts that Lemma \ref{lem6}. If $G\cong \mathrm{M}_{2}(2,1,1)$, then $r(G)>2$, this contradicts that Lemma \ref{lem3}. If $G\cong \mathrm{M}_{2}(2,2),~G\cong \mathrm{D}_{8}\times\mathrm{C}_{2}~ \mbox{or}~G\cong \mathrm{Q}_{8}\times\mathrm{C}_{2}$, then $Z(G)\cong \mathrm{C}_{2}\times \mathrm{C}_{2}$. This contradicts that Lemma \ref{cc}. If $G\cong \mathrm{D}_{8}*\mathrm{C}_{4}$, by calculation, $\delta_{\mathcal{CD}}(G)=18$. This contradicts the hypothesis. By calculation,  $\delta_{\mathcal{CD}}(\mathrm{M}_{2}(3,1))=6$.
%
%Case $3$: $2^{5}\leqslant|G|\leqslant2^{7}$.
Using Magma \cite{HUB} to check the SmallGroups Database \cite{WB}, we find the group $G$ with $\delta_{\mathcal{CD}}(G)\leqslant 6$ is $\mathrm{Q}_{8}$, $\mathrm{D}_{8}$ or $\mathrm{M}_{2}(3,1)$. By calculation,  $\delta_{\mathcal{CD}}(\mathrm{Q}_{8})=1$, $\delta_{\mathcal{CD}}(\mathrm{D}_{8})=5$ and $\delta_{\mathcal{CD}}(\mathrm{M}_{2}(3,1))=6$.

%there is no group that satisfies conditions.

\bigskip

\noindent {\bf Proof of Main Theorem.} The deduction can be made based on Lemmas \ref{lem1}, \ref{zhu2} and \ref{zhu3}.

%\begin{thm}
%Let $G$ be a ablian $2$-group. If $\delta_{\mathcal{CD}}(G)\leqslant 6$, then one of the following holds:
%
%$(1)$ $G\cong \mathrm{C}_{2^{k}},k=1,2,\ldots,6$;
%
%$(2)$ $G\cong \mathrm{C}_{2}\times\mathrm{C}_{2}$.
%\end{thm}
%
%
%
%\begin{thm}
%Let $G$ be a ablian $p$-group, where $p>2$. If $\delta_{\mathcal{CD}}(G)\leqslant p^{2}+p$, then one of the following holds:
%
%$(1)$ $G\cong \mathrm{C}_{p^{k}},k=1,2,\ldots,p^{2}+p$;
%
%$(2)$ $G\cong \mathrm{C}_{p}\times\mathrm{C}_{p}$;
%
%$(3)$ $G\cong \mathrm{C}_{p^2}\times\mathrm{C}_{p}$.
%\end{thm}
%
%
%\begin{thm}
%Let $G$ be a non-ablian $p$-group, where $p>2$. If $\delta_{\mathcal{CD}}(G)\leqslant p^{2}+p$, then one of the following holds:
%
%$(1)$ $G\cong \mathrm{M}_{p}(2,1)$;
%
%
%\end{thm}
%
%
%\begin{thm}
%Let $G$ be a non-ablian $2$-group. If $\delta_{\mathcal{CD}}(G)\leqslant 6$, then one of the following holds:
%
%$(1)$ $G\cong \mathrm{D}_{8}$;
%
%$(2)$ $G\cong \mathrm{Q}_{8}$;
%
%$(3)$ $G\cong \mathrm{M}_{2}(3,1)$;
%
%
%\end{thm}

\bigskip

%{\bf Acknowledgments}  I cordially thank the referee for detail
%reading and helpful comments, which helped me to improve the whole paper considerably.

\noindent{\bf Appendix.}

f:=function(p,n);

${\rm
P:=SmallGroupProcess(p^n);}$

X:=[];

repeat G:=Current(P);

if not IsAbelian(G) then

M:=Subgroups(G);

${\rm
T:=[x^{\backprime}length:x~ in~ M|\#(Centralizer(G,x^{\backprime}subgroup))*\#(x^{\backprime}subgroup)~ eq~}$\\
\vspace{0.5cm}
\hspace{1.5cm}${\rm \#G*\#(Center(G))];}$

\vspace{-0.5cm}${\rm S:=\&+T;}$

${\rm L:=[x^{\backprime}length:x~ in~ M];}$

${\rm K:=\&+L;}$

${\rm if~ K-S~ le~ 6~ then}$

${\rm
\_,a:=CurrentLabel(P);}$

${\rm
Append(^{\sim}X,a);}$

end if;

end if;

${\rm
Advance(^{\sim}P);}$

until IsEmpty(P);

return X;

end function;

time X:=f(2,3);

X;

${\rm
\#X;
}$

time X:=f(2,4);

X;

${\rm
\#X;
}$

time X:=f(2,5);

X;

${\rm
\#X;
}$

time X:=f(2,6);

X;

${\rm
\#X;
}$

time X:=f(2,7);

X;

${\rm
\#X;
}$

\end{document}